\documentclass[12pt,reqno]{amsart}
\usepackage{stmaryrd}
\usepackage{amscd}
\usepackage{pifont}
\usepackage{amsmath}
\usepackage{amsthm}
\usepackage{amsfonts}
\usepackage{amssymb}
\usepackage{latexsym}
\usepackage{amstext}
\usepackage{array}
\usepackage{graphicx}
\usepackage{slashbox}
\usepackage{url}

\date{}
\allowdisplaybreaks[4] \footskip=15pt
\renewcommand{\uppercasenonmath}[1]{}

\newtheorem{thm}[subsection]{Theorem}
\newtheorem{cor}[subsection]{Corollary}
\newtheorem{Def}[subsection]{Definition}
\newtheorem{lem}[subsection]{Lemma}
\newtheorem{remark}{Remark}

\newtheorem{prop}[subsection]{Proposition}
\newtheorem{exm}[subsection]{Example}
  
\newcommand{\bthm}{\begin{thm} }
\newcommand{\ethm}{\end{thm} }
\newcommand{\bpro}{\begin{prop}}
\newcommand{\epro}{\end{prop}}
\newcommand{\bdf}{\begin{Def}}
\newcommand{\edf}{\end{Def}}
\newcommand{\bexm}{\begin{exm}}
\newcommand{\eexm}{\end{exm}}
\newcommand{\blem}{\begin{lem}}
\newcommand{\elem}{\end{lem}}
\newcommand{\bpf}{\begin{proof}}
\newcommand{\epf}{\end{proof}}
\newcommand{\bcor}{\begin{cor}}
\newcommand{\ecor}{\end{cor}}

\newcommand{\bea}{\begin{eqnarray}}
\newcommand{\eea}{\end{eqnarray}}
\newcommand{\brem}{\begin{remark}}
\newcommand{\erem}{\end{remark}}

\input txdtools
\begin{document}

\begin{center}

{\huge  \bf On the periodicity of some Farhi  arithmetical
functions} \footnote{The first author is partially supported by the
 Grant No. 10571080 and 10871088 from NNSF of China.
 The second author  is partially supported by the Grant
 No. 10171046 from NNSF of China and Jiangsu
planned projects for postdoctoral research funds.}
\\

 \vskip 0.8cm
 {\small  Qing-Zhong Ji$^1$} and {\small  Chun-Gang Ji$^2$} \\

\end{center}

\vskip 0.8cm

 {\bf Abstract:}\quad
Let $k\in\mathbb{N}$. Let $f(x)\in \Bbb{Z}[x]$ be any polynomial
such that $f(x)$ and $f(x+1)f(x+2)\cdots f(x+k)$ are coprime in
$\mathbb{Q}[x]$. We call
$$g_{k,f}(n):=\frac{|f(n)f(n+1)\cdots f(n+k)|}{
\text{lcm}(f(n),f(n+1),\cdots,f(n+k))}$$ a Farhi arithmetic
function. In this paper, we  prove that $g_{k,f}$ is periodic.  This
generalizes the previous results of Farhi and Kane, and Hong and
Yang.

\vskip 0.8cm

\section{\bf Introduction}

Throughout this paper, let $\Bbb Q$, $\Bbb Z$ and $\Bbb N$ denote
the field of rational numbers, the ring of rational integers and the
set of nonnegative integers. Let ${\Bbb N}^*={\Bbb N}\setminus \{
0\}$. As usual, let $v_p$ denote the normalized $p$-adic valuation
of $\Bbb Q$, i.e., $v_p(a)=b$ if $p^b||a$.

It is known that an equivalent variation of the Prime Number Theorem
states that $\log {\rm lcm}(1,2,\cdots,n)\sim n$ as $n$ tends to
infinity (see e.g., \cite{HW}). One thus expects  that a better
understanding of the function ${\rm lcm}(1,2,\cdots,n)$ may entail a
deeper understanding of the distribution of the prime numbers. Some
progress has been made towards this direction.  Before we state our
main theorems, let us first give a short account on the recent
results in this subject.

In  his  pioneered paper \cite{Fa2},  Farhi
 introduced the following  arithmetic
functions
$$g_k(n):=\frac{n(n+1)\cdots(n+k)}{\text{
lcm} (n,n+1,\cdots,n+k)},\;\;  n\in {\Bbb N}^*.$$

Farhi proved that the sequence $(g_k)_{k\in {\Bbb N}}$ satisfies the
recursive relation:
\begin{eqnarray}g_k(n)=\gcd(k!, (n+k)g_{k-1}(n)),\;\;\forall n\in {\Bbb
N}^*.\label{001}\end{eqnarray} Using this relation, Farhi proved
\bthm{\rm(\cite{Fa2})} The function $g_k\;(k\in {\Bbb N})$ is
periodic and $k!$ is a period of $g_k.$ \ethm
 An interesting question is  how to determine the least
period of $g_{k}$ (see \cite{Fa2}).  In  \cite{HY}, by using
(\ref{001}) and $g_k(1)|g_k(n)$ for any positive integer $n$, Hong
and Yang gave a partial answer to this question. A complete solution
to the question was given by Farhi and Kane in \cite{FK}. They
proved \bthm{\rm(\cite{FK}, Theorem 3.2)}\label{100} The least
period $T_k$ of $g_k$ is given by
\begin{eqnarray}T_{k}=\prod_{p \; prime,\; p\leq
k}p^{\delta_p(k)},\label{10-1}\end{eqnarray} where
$$\delta_p(k)=\left\{ \begin{array}{ll}0, & {\rm if } \;v_p(k+1)\geq
\max\limits_{1\leq i\leq k}\{v_p(i)\},\\
\max\limits_{1\leq i\leq k} \{v_p(i)\}, & {\rm
otherwise.}\end{array}\right.$$ \ethm

Let $g(n)$ be an arithmetic function defined on the set
$\mathbb{Z}\backslash A,$ where $A$ is a finite subset of
$\mathbb{Z}.$ If there exists an integer $T$ such that $g(n)=g(n+T)$
for all $n,\;n+T\in\mathbb{Z}\backslash A,$  then it is clear that
the arithmetic function $g(n)$ can be extended a periodic function
defined on  all the integers $\mathbb{Z}.$

Throughout this paper, let $f(x)\in \mathbb{Z}[x]$ be a polynomial
such that $$\gcd(f(x), f(x+1)f(x+2)\cdots f(x+k))=1$$ in
$\mathbb{Q}[x]$. Let $k$ be a nonnegative integer.  Denoted by
\begin{eqnarray}Z_{k, f}:=\{n\in \mathbb{Z}\ | \ f(n+i)=0 \mbox{ for
some } 0\leq i\leq k\}.\label{002-1}\end{eqnarray} Then $Z_{k, f}$
is a finite subset of $\mathbb{Z}$. Set
\begin{eqnarray}g_{k, f}(n)=\frac{|f(n)f(n+1)\cdots f(n+k)|}{{\rm
lcm}(f(n),f(n+1),\cdots,f(n+k))},\label{002}\end{eqnarray} for $n\in
\mathbb{Z}\backslash Z_{k, f}$. We call $g_{k,f}(n)$ a Farhi
arithmetic function. In $\S 3$, we will prove

\bthm\label{3} Let $k$ be a nonnegative integer and $f(x)\in
\mathbb{Z}[x]$ be a polynomial such that $\gcd(f(x),
f(x+1)f(x+2)\cdots f(x+k))=1$ in $\mathbb{Q}[x].$ Then the
arithmetic function $g_{k,f}$ can be extended to a periodic
arithmetic function defined on all the integers$.$ \ethm

By assumption of $f(x)$ in Theorem \ref{3}, for any $1\leq i \leq
k$, there exist polynomials $a_i(x),b_i(x)\in\mathbb{Z}[x]$ and the
smallest positive integer $C_i$ such that
$$a_i(x)f(x)+b_i(x)f(x+i)=C_i.$$
Let $C$ be the least common multiple of the $C_i$'s, i.e.,
$$C={\rm lcm}(C_1, C_2, \cdots, C_k).$$
In the proof of Theorem \ref{3}, we will prove

\bthm\label{4}  Let $T_{k, f}$ denote the  least period of
$g_{k,f}.$ Then $T_{k, f}|C.$ \ethm

Let $p$ be a prime. Define the arithmetic function
 $h_{k,f,p}$ by
 \begin{eqnarray}h_{k,f,p}(n):=v_p(g_{k,f}(n)).\label{0010-1}\end{eqnarray}
If $p\nmid C$, using the definition of $g_{k, f}$, then we have
$h_{k,f,p}(n)=1$ for any $n\in \Bbb Z$. If $p|C$, then $h_{k,f,p}$
is a periodic function by  theorem 1.3. Set
\begin{eqnarray}S_n:=\{n,n+1,\cdots,n+k\},\;n\in\mathbb{Z}\label{0010-2}\end{eqnarray}
and
\begin{eqnarray}e_p:=\max\{v_p(\gcd(|f(n)|,|f(n+i)|))\;|\;1\leq n\leq p^{v_p(C)}, 1\leq i\leq
k\}.\label{0010}\end{eqnarray}

In $\S 4$, we will prove

\bthm\label{5} For any prime $p$, let $T_{k, f, p}$ be the least
period of the arithmetic function $h_{k,f,p}.$ Then

{\rm (i)} $p^{e_p}$ is a period of $h_{k,f,p}$ and $T_{k, f,
p}|p^{e_p}.$

{\rm (ii)} $T_{k,f,p}=1$ if and only if for any $1\leq n\leq
p^{e_p},$ we have $$v_p(\gcd(|f(n)|,|f(n+k+1)|))\geq \underset{1\leq
i\leq k}\max\{v_p(f(n+i))\}$$ or
$$v_p(f(n))=v_p(f(n+k+1))<\underset{1\leq i\leq k}\max\{v_p(f(n+i))\}.$$

{\rm (iii)} Let $1\leq e\leq e_p.$ Suppose that $p^e$ is a period of
$h_{k,f,p}.$ Then $T_{k,f,p}=p^e$ if and only if there exists an
integer $n_0:\;1\leq n_0\leq p^{e}$ such that
 the following inequality holds$:$
$$\sum_{t=e}^{e_p}\max\{0, \#
\{m\in S_{n_0}\; |\; p^t|f(m)\}-1\}$$
$$\not=\sum_{t=e}^{e_p}\max\{ 0, \# \{m\in S_{n_0}\; |\;
p^t|f(m+p^{e-1})\}-1\}.$$ In particular, $T_{k,f,p}=p^{e_p}$ if and
only if  there exists an integer $n_0:\;1\leq n_0\leq p^{e_p}$ such
that the following inequality holds$:$
 $$\#\{m\in
S_{n_0}\;|\;p^{e_p}|f(m)\}\not=\#\{m\in S_{n_0}\;|\;
p^{e_p}|f(m+p^{e_p-1})\}.$$ \ethm

{\bf Remark: } Let $T_{k, f, p}$ be  the least period of $h_{k, f,
p}$ for any prime $p.$ Then
$$T_{k, f}=\prod_{p}T_{k, f, p}$$
(This infinite product is meaningful, for almost all its terms are
equal to 1).

As an application of Theorem \ref{5}, in $\S 5$ we will give a new
different proof of Theorem 3.2 of \cite{FK}.

\bcor\label{17}\ \ Let $k\in\mathbb{N}$ and $f(x)=x.$ Then the least
 period $T_{k,f}$ of the {\rm Farhi} arithmetic function $g_{k,f}$ is
given by the  formula {\rm (\ref{10-1}).} \ecor

 \vskip 2mm

Let $a,b\in\mathbb{Z}$ be any integer such that ${\rm gcd}(a,b)=1$
and $a>0.$ Let $f(x)=ax+b.$ By Theorem \ref{3}, we know  that the
Farhi arithmetic function $$g_{k,
ax+b}(n)=\frac{|(an+b)(a(n+1)+b)\cdots (a(n+k)+b)|}{{\rm
lcm}(an+b,a(n+1)+b,\cdots,a(n+k)+b)}$$ can be extended to a periodic
arithmetic function defined on all the integers. Now we define the
arithmetical function $g_{k,a}$ by
$$g_{k,a}(n)=\frac{|n(n+a)\cdots(n+ka)|}{\text{lcm}(n,n+a,\cdots,n+ka)}.$$
When $a=1$, the arithmetical function $g_{k,1}$ is the arithmetical
function $g_k$ defined by Farhi. It is clear that
\begin{eqnarray}g_{k,ax+b}(n)=g_{k,a}(na+b).\label{10001}\end{eqnarray} Hence the function $g_{k,a}$ also
can be extended to a periodic arithmetic function defined on all the
integers. In $\S 6,$ we shall prove the following results:

\bthm\label{18} Let $a$, $k$ be any two positive integers$.$ Then
the following assertions hold$.$

{\rm(i)} The positive integer $a\cdot \text{lcm}(1, 2, \cdots, k)$
is a period of $g_{k, a}.$

{\rm(ii)} A positive integer $S$ is a period of $g_{k,a}$ if and
only if $S=aT,$ where $T$ is a period of $g_k;$

{\rm(iii)} Consequently, the least period $T_k(a)$ of $g_{k,a}$ is
$aT_k(1)=aT_k,$ where $T_k(1)=T_k$ is the least period of $g_k.$
\ethm

\vskip 2mm

By (\ref{10001}) and Theorem \ref{18}, we have the following result:

\bcor \label{19}Let $a$, $k$ be any two positive integers and let
$b\in\mathbb{Z}$ be any integer such that ${\rm gcd}(a,b)=1.$ Then
the least period $T_{k,ax+b}$ of the Farhi arithmetic function
$g_{k,ax+b}$ is given by the following formula:
\begin{eqnarray}T_{k,ax+b}=\prod_{p \; prime,\; p\leq
k}p^{\delta_p(k)},\label{10-1}\end{eqnarray} where
$$\delta_p(k)=\left\{ \begin{array}{ll}0, & {\rm if } \;v_p(k+1)\geq
\max\limits_{1\leq i\leq k}\{v_p(i)\}\;{\rm or}\;p|a,\\
\max\limits_{1\leq i\leq k} \{v_p(i)\}, & {\rm
otherwise.}\end{array}\right.$$

 \ecor

 In $\S 7$, we will give some examples.

\section{\bf Two Basic Lemmas }

\begin{lem}
\label{1} Let $a_1, a_2, \cdots, a_n$ and $b_1, b_2, \cdots, b_n$ be
any $2n$ positive integers$.$ If
 $\gcd(a_i, a_j)=\gcd(b_i, b_j)$ for any $1\leq i<j\leq
 n,$
then \begin{eqnarray}\label{003}\frac{a_1a_2\cdots a_n}{{\rm
lcm}(a_1,a_2,\cdots,a_n)} =\frac{b_1b_2\cdots b_n}{{\rm
lcm}(b_1,b_2,\cdots,b_n)}.\end{eqnarray}
\end{lem}

\begin{proof}
 Let $p$ be any prime. It suffices to show that
the following equality holds:
\begin{eqnarray}\label{004}v_p\left(\frac{a_1a_2\cdots
a_n}{{\rm lcm}(a_1,a_2,\cdots,a_n)}\right)
=v_p\left(\frac{b_1b_2\cdots b_n}{{\rm
lcm}(b_1,b_2,\cdots,b_n)}\right).
\end{eqnarray}
i.e.,
\begin{eqnarray}\sum_{i=1}^nv_p(a_i)-\max_{1\leq i\leq n}\{v_p(a_i)\}
= \sum_{i=1}^nv_p(b_i)-\max_{1\leq i\leq
n}\{v_p(b_i)\}.\label{a0}\end{eqnarray} By the assumption of $a_i$'s
and $b_i$'s, it suffices to show that
\begin{eqnarray}\sum_{i=1}^nv_p(a_i)-\max_{1\leq i\leq n}\{v_p(a_i)\}
\leq \sum_{i=1}^nv_p(b_i)-\max_{1\leq i\leq
n}\{v_p(b_i)\}.\label{0001}\end{eqnarray} Without loss of
generality, we assume that $v_p(a_1)\leq v_p(a_2)\leq \cdots \leq
v_p(a_{n-1})\leq v_p(a_n).$ Then for any $1\leq i\leq n-1$, we have
$$v_p(a_i)=v_p({\rm gcd}(a_i, a_n))=v_p({\rm gcd}(b_i, b_n))
\leq \min\{v_p(b_i), v_p(b_n)\}.$$ Hence for any $1\leq i\leq n-1$,
we have $v_p(a_i)\leq v_p(b_i)$ and $v_p(a_i)\leq v_p(b_n)$. Let
$v_p(b_k)=\max_{1\leq i\leq n}\{v_p(b_i)\}$. Then
$$\sum_{i=1}^{n-1}v_p(a_i)\leq v_p(b_1)+\cdots+v_p(b_{k-1})+v_p(b_n)+
v_p(b_{k+1})+\cdots+v_p(b_{n-1}).$$ So $(\ref{0001})$ is true. This
completes the proof of Lemma 2.1. \qed
\end{proof}

\begin{lem}
\label{2} Let $k$ be a positive integer and $f(x)\in \mathbb{Z}[x]$
be any polynomial such that $$\gcd (f(x), f(x+1)f(x+2)\cdots
f(x+k))=1$$ in $\mathbb{Q}[x].$ Then $d_{i}(n)=\gcd (|f(n)|,
|f(n+i)|)$ is periodic for any $1\leq i\leq k.$
\end{lem}

\begin{proof}
By assumption, for any $1\leq i\leq k$, $f(x)$ and $f(x+i)$ are
coprime in $\mathbb{Q}[x]$, hence there exist $a_i(x),
b_i(x)\in\mathbb{Z}[x]$ and the smallest positive integer $C_i$ such
that
\begin{eqnarray}\label{006}a_i(x)f(x)+b_i(x)f(x+i)=C_i.\end{eqnarray}
Hence for all $m\in\mathbb{Z},$ we have
\begin{eqnarray}\label{2004}a_i(m)f(m)+b_i(m)f(m+i)=C_i.\end{eqnarray}
In the following, we will prove
\begin{eqnarray}\label{2001}d_i(n)=d_i(n+C_i), \quad n\in \mathbb{Z}.\end{eqnarray}
Let $d_i=d_i(n)$ and $d_i'=d_i(n+C_i)$. Then $d_i|f(n)$ and
$d_i|f(n+i)$. By (\ref{2004}), we have $d_i|C_i.$ Hence
$d_i|f(n+C_i)$ and $d_i|f(n+i+C_i)$. Therefore $d_i|d_i'$.
Similarly, we have $d_i'|d_i$. Hence $d_i=d_i'$, i.e., (\ref{2001})
is true. This completes the proof of Lemma 2.2. \qed
\end{proof}

\section{\bf The Proofs of Theorem \ref{3} and Theorem \ref{4}}

{\bf Proof of Theorem \ref{3}.} By the definition (\ref{002-1}),
$Z_{k,f}$ is a finite set and $g_{k,f}$ is well defined on the set
$\mathbb{Z}\setminus Z_{k,f}.$  First we prove that $g_{k,f}$ is
periodic on the set $\mathbb{Z}\setminus Z_{k,f}.$ For $1\leq i\leq
k,$ by Lemma \ref{2}, $d_i(n)=\gcd(|f(n)|,|f(n+i)|)$ is periodic.
Let $T_i$ be the least period of $d_i$. Then
$$d_i(n)=d_i(n+T_i), \quad \mbox{ for any } n\in \mathbb{Z}.$$
Hence by the proof of Lemma \ref{2}, we have that $T_i|C_i,$ where
$C_i$ is defined by (\ref{006}). Denote by $T$ (resp. $C$) the least
common multiple of the $T_i$'s (resp. $C_i$'s), $i=1,2,\cdots,k$.
Then $T|C$ and for any $1\leq i\leq k$, we have
$$d_i(n)=d_i(n+T)\quad \mbox{ for } n\in \mathbb{Z}.$$
Hence for any $0\leq i<j\leq k$, we have
$$d_{j-i}(n+i)=d_{j-i}(n+i+T) \quad \mbox{ for } n\in \mathbb{Z},$$
that is
$$\gcd(|f(n+i)|,|f(n+j)|)=\gcd(|f(n+i+T)|,|f(n+j+T)|).$$
So by Lemma \ref{1} and the definition of $g_{k, f}$, we obtain that
 $g_{k,f}(n)=g_{k,f}(n+T)$ for any $n$ and $n+T \in \mathbb{Z}\backslash Z_{k, f}$.
 Hence $g_{k,f}(n)$ is periodic and $T$ is a period of $g_{k,f}$.

If $n\in Z_{k, f}$, then there exist a positive integer $a$ such
that $n+aT\not\in Z_{k,f}$.  Hence the function $g_{k,f}$ can be
extended to $g_{k,f}:\;\mathbb{Z}\longrightarrow \mathbb{Z}$ defined
at $n\in Z_{k,f},$ by
$$g_{k,f}(n)=g_{k,f}(n+aT).$$
This completes the proof of Theorem 1.3. \qed

\vskip 2mm

{\bf Proof of Theorem \ref{4}.} It is obvious that the property
$T|C$ implies that $C$ is a multiple of the least period $T_{k,f}$
of $g_{k,f}$. this completes the proof of Theorem 1.4. \qed

\section{\bf The Proof of Theorem \ref{5}}

Notations as previous sections.

\bpf \quad (i) By the definitions of $h_{k, f, p}$ and $g_{k, f}$,
it suffices to show that $v_p(g_{k, f}(n))=v_p(g_{k, f}(n+p^{e_p}))$
for any $n\in \mathbb{Z}\setminus Z_{k,f},$ i.e.,
$$\sum_{i=0}^{k}v_p(f(n+i))-\max_{0\leq i\leq k}\{v_p(f(n+i))\}$$
$$=\sum_{i=0}^{k}v_p(f(n+i+p^{e_p}))-\max_{0\leq i\leq
k}\{v_p(f(n+i+p^{e_p}))\}.$$ Let
$$e_{ij}=v_p(\gcd(|f(n+i)|, |f(n+j)|))$$ and
$$e_{ij}'=v_p(\gcd(|f(n+i+p^{e_p})|, |f(n+j+p^{e_p})|))$$ for any
$0\leq i<j\leq k$. By the proof
of Lemma 2.1, it suffices to show that
\begin{eqnarray}e_{ij}=e_{ij}'.
\label{2003}\end{eqnarray} By the assumption of $f(x)$, we have
$$a_{j-i}(m)f(m)+b_{j-i}(m)f(m+j-i)=C_{j-i},\quad m\in \mathbb{Z}$$
Let $m=n+i$. We have $p^{e_{ij}}|f(n+i)$ and $p^{e_{ij}}|f(n+j)$, so
$p^{e_{ij}}|C_{j-i}$. Hence $e_{ij}\leq e_p$ by the definition of
$e_p$. So $p^{e_{ij}}|f(n+i+p^{e_p})$, $p^{e_{ij}}|f(n+j+p^{e_p})$.
Therefore $e_{ij}\leq e_{ij}'$. Similarly, we have $e_{ij}'\leq
e_{ij}$. Hence (\ref{2003}) is true. It is easy to see that $T_{k,
f, p}|p^{e_p}$.

\vskip5mm (ii) By (i) of Theorem 1.5, we know that $h_{k, f, p}$ is
periodic and $p^{e_p}$ is a period. So $T_{k, f, p}=1$ if and only
if $h_{k, f, p}(n)=h_{k, f, p}(n+1)$ for any $1\leq n\leq p^{e_p}$.
By the definition of $g_{k, f}$, we have $T_{k, f, p}=1$ if and only
if for any $1\leq n\leq p^{e_p}$,
$$\sum_{i=0}^{k}v_p(f(n+i))-\max_{0\leq i\leq k}\{v_p(f(n+i))\}$$
$$=\sum_{i=1}^{k+1}v_p(f(n+i))-\max_{1\leq i\leq
k+1}\{v_p(f(n+i))\}. $$ Hence $T_{k,f,p}=1$ if and only if for any
$1\leq n\leq p^{e_p}$, $$v_p(\gcd(|f(n)|,|f(n+k+1)|))\geq
\underset{1\leq i\leq k}\max\{v_p(f(n+i))\}$$ or
$$v_p(f(n))=v_p(f(n+k+1))<\underset{1\leq i\leq k}\max\{v_p(f(n+i))\}.$$

\vskip5mm (iii) Let $1\leq e\leq e_p.$ Suppose that $p^{e}$ is a
period of $h_{k,f,p}.$ Hence $p^{e}$ is the least period of
$h_{k,f,p}$ if and only if $p^{e-1}$ is not a period $h_{k,f,p}.$
Therefore $p^{e}$ is the least period of $h_{k,f,p}$ if and only if
there exists an integer $n_0:\;1\leq n_0\leq p^{e}$ such that the
following inequality holds:
$$h_{k,f,p}(n_0)\not=h_{k,f,p}(n_0+p^{e-1}).$$
 By definition (\ref{0010-1}), we have
\begin{eqnarray}\label{0012}
\begin{array}{ll}h_{k,f,p}(n_0)&=\sum\limits_{i=0}^kv_p(f(n_0+i))-\underset{0\leq i\leq
k}{\max} \{v_p(f(n_0+i))\}\\[4mm]&=\sum\limits_{t=1}^{\infty}\max \{0,\; \# \{m\in
S_{n_0}\;|\;p^t|f(m)\}-1\},\end{array}
\end{eqnarray} and
\begin{eqnarray}\label{0013}
h_{k,f,p}(n_0+p^{e-1})=\sum_{t=1}^{\infty}\max \{0,\; \# \{m\in
S_{n_0}\;|\;p^t|f(m+p^{e-1})\}-1\}.
\end{eqnarray}
Remark: This infinite sum is meaningful, for almost all its terms
are equal to 0, and $\#\emptyset=0$. On the other hand, when $t\leq
e-1,$ we know that $p^t|f(m)$ if and only if $p^t|f(m+p^{e-1}).$
Hence by the definition of $e_p$, the inequality
$h_{k,f,p}(n_0)\not=h_{k,f,p}(n_0+p^{e-1})$ holds if and only if the
inequality $$\sum_{t=e}^{e_p}\max\{0, \# \{m\in S_{n_0}\; |\;
p^t|f(m)\}-1\}$$
$$\not=\sum_{t=e}^{e_p}\max\{ 0, \# \{m\in S_{n_0}\; |\;
p^t|f(m+p^{e-1})\}-1\}$$ holds. In particular, $T_{k,f,p}=p^{e_p}$
if and only if  there exists an integer $n_0:\;1\leq n_0\leq
p^{e_p}$ such that the following inequality holds:
 $$\#\{m\in
S_{n_0}\;|\;p^{e_p}|f(m)\}\not=\#\{m\in S_{n_0}\;|\;
p^{e_p}|f(m+p^{e_p-1})\}.$$ This completes the proof of Theorem
\ref{5}. \qed\epf

\section{\bf The proof of Corollary \ref{17}}

\begin{proof}
 When $k=0,$ then $g_{0,f}=1.$ Let $k\geq 1.$ For $1\leq i
\leq k,$ we have $C_i=i$ and $C={\rm lcm}(1,2,\cdots,k).$ Hence we
obtain
$$T_{k,f}=\prod_{p \; prime, p\leq k}T_{k,f,p}.$$
Let $p\leq k$ be a prime, it suffices to prove the following
statements:

(I) $T_{k,f,p}=1$ if $v_p(k+1)\geq \underset{1\leq i\leq k}{\rm
max}\{v_p(i)\}.$

(II) $T_{k,f,p}=p^{v_p(C)}$ if $v_p(k+1)< \max\limits_{1\leq i\leq
k}\{v_p(i)\}.$

We first prove (I). As $e_p=v_p(C)=\max\limits_{1\leq i\leq
k}\{v_p(i)\}$, by assumption $v_p(k+1)\geq e_p,$ we have
$v_p(k+1)=e_p$ or $e_p+1$.

Case (a), $1\leq n\leq p^{e_p}-1$, then $e=v_p(n)<e_p$. Hence
$v_p(n)=v_p(n+k+1)$ and $n=p^{e}n_1$, $p\nmid n_1.$ Set
$i=p^{e}i_0$, $1\leq i_0\leq p-1$ such that $p|(n_1+i_0)$. Then
$1\leq i\leq k$ and $v_p(n+i)>v_p(n)$. Hence
\begin{eqnarray}\label{0014}v_p(n)=v_p(n+k+1)<\max\limits_{1\leq i\leq
k}\{v_p(n+i)\}.\end{eqnarray}

Case (b1), $n=p^{e_p}$ and $v_p(k+1)=e_p+1$. We  have
$k+1=p^{e_p+1}$. Let $i=p^{e_p}(p-1)$. Then $1\leq i\leq k$ and
$v_p(n+i)=e_p+1>e_p$. Hence
\begin{eqnarray}\label{0015}v_p(n)=v_p(n+k+1)<\max\limits_{1\leq
i\leq k}\{v_p(n+i)\}.\end{eqnarray}

Case (b2), $n=p^{e_p}$ and $v_p(k+1)=e_p$. We  have $k+1=p^{e_p}u$,
where $2\leq u\leq p-1$. Hence $k=up^{e_p}-1$. If $i>0$ and
$v_p(n+i)>e_p$, then $i\geq p^{e_p}(p-1)>k$. Hence
$\max\limits_{1\leq i\leq k}\{v_p(n+i)\}\leq e_p.$  Therefore
\begin{eqnarray}\label{0016}v_p(\gcd(|n|, |n+k+1|))\geq \max\limits_{1\leq i\leq
k}\{v_p(n+i)\}.\end{eqnarray} By (\ref{0014}), (\ref{0015}),
(\ref{0016}) and using (ii) of Theorem 1.5, we have $T_{k,f,p}=1$.

(II) Note that $e_p=v_p(C).$ Hence
$$k=a_0+a_1p+\cdots+a_{e_p}p^{e_p},\;\;0\leq a_i\leq p-1,\;\;i=0,1,\cdots, e_p,\;a_{e_p}\not=0.$$
It is easy to show  that the inequality $v_p(k+1)\geq e_p$ holds if
and only if $a_0=a_1=\cdots=a_{e_p-1}=p-1.$

Assume that the inequality $v_p(k+1)<e_p=v_p(C)=\underset{1\leq
i\leq k}{\rm max}\{v_p(i)\}$ holds. Then there exists an integer
$r:\;0\leq r\leq e_p-1$ such that the following conditions hold:
$$0\leq a_r\leq p-2\;{\rm and}\; a_{r+1}=\cdots=a_{e_p-1}=p-1.$$ Set
\[
n_0=\left\{
\begin{array}{ll}
p^{e_p},&{\rm if}\;r=e_p-1;\\
(p-1-a_r)p^{r},&{\rm if}\;0\leq r\leq e_p-2.
\end{array}
\right.\] Then we have
$$\# \{m\in
S_{n_0}\;\;|\;\;p^{e_p}|m\}=\left\{
\begin{array}{ll}
a_{e_p}+1,&{\rm if}\;r=e_p-1,\\
a_{e_p},&{\rm if}\;0\leq r\leq e_p-2;\end{array}\right.$$and
$$\# \{m\in
S_{n_0}\;\;|\;\;p^{e_p}|(m+p^{e_p-1})\}=\left\{
\begin{array}{ll}a_{e_p},&{\rm if}\; r=e_p-1,\\ a_{e_p}+1,&{\rm if}\;0\leq r\leq e_p-2.\end{array}\right.$$
By the (iii) of Theorem \ref{5}, we know that $p^{e_p}=p^{v_p(C)}$
is the least period of $h_{k,f,p}.$ This completes the proof of
Corollary \ref{17}.\qed

\end{proof}

\vskip 2mm

\section{\bf The proof of Theorem \ref{18}}

\bpf$\;$ (i) Set $S=a\cdot \text{lcm}(1, 2, \cdots, k)$. Let $n$ be
any positive integer. For any $0\leq i<j\leq k$, it is clear that
$\text{gcd}(n+ia, n+ja)=\text{gcd}(n+S+ia, n+S+ja)$. Hence
$g_{k,a}(n+S)=g_{k,a}(n)$ follows from Lemma \ref{1}.

(ii)   Suppose $S$ is a period of $g_{k, a}$. Then $g_{k,
a}(n)=g_{k, a}(n+S)$ for all $n\in {\Bbb N}^*$. In particular, we
have $g_{k, a}(na)=g_{k,a}(na+S)$. Since
$$g_{k,a}(na)=\frac{na\cdot (na+a)\cdot \cdots (na+ka)}{\text{lcm}(na,
na+a, \cdots, na+ka)}=\frac{n(n+1)\cdots
(n+k)}{\text{lcm}(n,n+1,\cdots,n+k)}\cdot a^k$$
 and
 $$g_{k,a}(na+S)=\frac{(na+S)(na+a+S)\cdots(na+ka+S)}{\text{lcm}(na+S, na+a+S,\cdots,na+ka+S)},$$
we have \begin{eqnarray}g_k(n)\cdot
a^k=g_{k,a}(na)=g_{k,a}(na+S)\label{7003}\end{eqnarray} and
\begin{eqnarray}g_k(n)\cdot a^k
=\frac{(na+S)(na+a+S)\cdots(na+ka+S)}{\text{lcm}(na+S,
na+a+S,\cdots,na+ka+S)}.\label{7004}\end{eqnarray} We claim that
$a|S$. Let $\text{gcd}(a, S)=d$, $a=a_1d$, $S=S_1d$. Then
$\text{gcd}(a_1, S_1)=1$. By using (\ref{7004}), we have
$$ a_1^k \, |\, (na_1+S_1)(na_1+a_1+S_1)\cdots(na_1+ka_1+S_1).$$
Because $\text{gcd}(a_1, na_1+ia_1+S_1)=1$ for any $0\leq i\leq k$,
we have $a_1=1$. Hence $a|S$. Let $S=aT$. Then using (\ref{7003}),
we have
$$g_k(n)\cdot a^k=g_{k,a}(na)=g_{k,a}(na+aT)=g_k(n+T)\cdot a^k.$$
Hence $g_k(n+T)=g_k(n)$ for all $n\in {\Bbb N}^*$, i.e., $T$ is a
period of $g_k(n)$.

Conversely, suppose $T$ is a period of $g_k(n).$ Let $n$ be any
positive integer. If $d=\text{gcd}(n,a)$, $n=n_1d$, $a=a_1d$, then
$\text{gcd}(n_1, a_1)=1$ and
$$g_{k,a}(n)=g_{k,a_1}(n_1)\cdot d^k,\quad g_{k,a}(n+aT)=g_{k,a_1}(n_1+a_1T)\cdot d^k.$$
Hence, without loss of generality, we assume that $(n,a)=1$.
Therefore \begin{eqnarray}\text{gcd}(a, g_{k,a}(n))=1,\quad
\text{gcd}(a, g_{k,a}(n+aT))=1.\label{7005}\end{eqnarray} Hence by
using (\ref{7005}), we have
$$g_{k, a}(n)=g_{k, a}(n+aT)$$
if and only if \begin{eqnarray}v_p(g_{k, a}(n))=v_p(g_{k, a}(n+aT)),
\label{7006}\end{eqnarray} for any prime $p\nmid a$.

Let $p$ be any prime such that $p\nmid a$ and $N$ be a positive
integer greater than $v_p(k!)$. Then there exists a unique positive
integer $m$ such that $1\leq m<p^N$ and \begin{eqnarray}ma\equiv 1 (
\text{ mod }p^N).\label{7007}\end{eqnarray} Let $p$, $n$, $a$, $m$
as above and $0\leq i<j\leq k$. then for any integer $l$, there are
\begin{eqnarray}v_p(\text{gcd}(n+al+ai, (j-i)a))= v_p(\text{gcd}(mn+l+i, j-i)).
\label{7008}\end{eqnarray} Let \begin{eqnarray}\text{gcd}(n+al+ai,
(j-i)a)=p^{x_{ij}}w, \quad p\nmid w \label{7009}\end{eqnarray} and
\begin{eqnarray}\text{gcd}(mn+l+i, j-i)=p^{y_{ij}}u, \quad p\nmid u.
\label{7010}\end{eqnarray} Then by (\ref{7009}), there exists $s_1$,
$t_1\in {\Bbb Z}$ such that $(n+al+ai)s_1+(j-i)at_1=p^{x_{ij}}w.$
Multiplied by $m$ on both sides, we have
$(mn+aml+ami)s_1+(j-i)mat_1=p^{x_{ij}}mw.$ Using (\ref{7007}), we
have $(mn+l+i)s_1+(j-i)t_1=p^{x_{ij}}mw-p^N\delta.$ By (\ref{7010}),
we have $y_{ij}\leq x_{ij}$. Conversely, by (\ref{7010}), there
exists $s_2$, $t_2\in {\Bbb Z}$ such that
$(mn+l+i)s_2+(j-i)t_2=p^{y_{ij}}u.$ Multiplied by $a$ on both sides,
we have $(mna+al+ai)s_2+(j-i)at_2=p^{y_{ij}}au.$ Similarly, we have
$x_{ij}\leq y_{ij}$. So $x_{ij}= y_{ij}$ and (\ref{7008}) is true.
Let $l=0$ and $T$. By using (\ref{7008}) we have
$$v_p(\text{gcd}(n+ai, n+aj))=v_p(\text{gcd}(mn+i, mn+j))$$
and
$$v_p(\text{gcd}(n+aT+ai, n+aT+aj))=v_p(\text{gcd}(mn+T+i, mn+T+j))$$
for any $0\leq i<j\leq k$. By the proof of Lemma \ref{1}, we have
$$v_p(g_{k,a}(n))=v_p(g_k(mn)),\quad
v_p(g_{k,a}(n+aT))=v_p(g_k(mn+T)).$$ So using $g_k(mn)=g_k(mn+T)$,
we have $v_p(g_{k,a}(n+aT))=v_p(g_{k,a}(n))$ for any prime $p$ such
that $p\nmid a$.  Hence $g_{k,a}(n+aT)=g_{k,a}(n)$ and $aT$ is a
period of $g_{k,a}(n)$.

 The proof of (iii) follows from (ii). This completes the proof of
Theorem \ref{18}.\epf\qed

\vskip 2mm

{\bf Proof of Corollary \ref{19}} (i)  Assume that $p|a$. Then it is
clear that the equality $v_p(g_{k,ax+b}(n))=0$ holds for any integer
$n$ when $g_{k,ax+b}(n)$ is well defined.

 (ii) Assume that $p$ is not a prime factor of $a.$ By the formula (\ref{10001}), we have that
$T_p$ is a period of $v_p(g_{k,ax+b}(n))$ if and only if $aT_p$ is a
period of $v_p(g_{k,a}(n))$. Hence, by Theorem \ref{18}, we have
that $T_p$ is a period of $v_p(g_{k,ax+b}(n))$ if and only if $T_p$
is a period of $v_p(g_{k}(n))$. Therefore Corollary \ref{19} is
obtained by Theorem \ref{100}.

\vskip 2mm

\section{\bf Examples}

\blem\label{11} Let $f_1(x)=f_2(x)^r,$ where $r\geq 1$ is an
integer$.$ Then $T_{k,f_1}=T_{k,f_2}.$ \elem

\bpf \quad By (\ref{002}), we have $g_{k,f_1}(n)=g_{k,f_2}(n)^r.$
Hence the result is obvious.\qed\epf

\vskip 2mm

{\bf Example 1.} Let $f(x)=x^r,\;r\geq 1.$ Then by Lemma \ref{11},
we have $T_{k,x^r}=T_{k,x},$ where $T_{k,x}$ is given by the formula
(\ref{10-1}).

\vskip2mm
 {\bf Example 2.} Let $f(x)=x^2+b.$ For $1\leq i\leq k,$ we
have
$$(2x+3i)(x^2+b)+(-2x+i)((x+i)^2+b)=i(i^2+4b),\;\;{\rm if}\;i\;\;{\rm is\;\;odd}.                 $$
$$(x+3j)(x^2+b)+(-x+j)((x+2j)^2+b)=4j(j^2+b),\;{\rm if}\;i=2j.$$
Hence
$$C_i=\left\{\begin{array}{ll}i(i^2+4b),&{\rm if}\;i\;{\rm is\;odd,}\\
4j(j^2+b),&{\rm if}\;i=2j.\end{array}\right.$$ Hence, given any
$k\in\mathbb{N}$ and $b\in\mathbb{Z}$, by Theorem \ref{5}, we can
determine the least period $T_{k,f}$ of the arithmetic function
$g_{k,f}.$ For $1\leq k\leq 6$ and $1\leq b\leq 6$, Table I gives
the $T_{k,f}$'s.

\vskip 1cm
\newsavebox{\tablebox}
\begin{lrbox}{\tablebox}
\begin{tabular}{|c|c|c|c|c|c|c|}
\multicolumn{7}{c}{\bf TABLE I}\\
\multicolumn{7}{c}{$\mbox{ The least period } T_{k,f} \mbox{ of } g_{k,f} \mbox{ with } f(x)=x^2+b$}\\

\multicolumn{7}{c}{}\\

\hline \backslashbox{$f(x)$}{$k$}&1&2&3&4&5&6\\
\hline $x^2+1$&5&2$\cdot$ 5&2$\cdot$ 3$\cdot$ 5$\cdot$ 13&2$\cdot$
3$\cdot$ 5$\cdot$ 13&2$\cdot$ 3$\cdot$ 5$\cdot$ 13$\cdot$
29&2$\cdot$
3$\cdot$ 5$\cdot$ 13$\cdot$ 29\\

\hline $x^2+2$&$3^2$&$2\cdot3^2$&$2\cdot 3^2\cdot 17$ &$2\cdot
3^2\cdot17$&$2\cdot 3^2\cdot 5\cdot 11\cdot 17$&$2\cdot 3^2\cdot 5\cdot 11\cdot 17$\\

\hline $x^2+3$&$13$&$2\cdot13$&$2\cdot 3\cdot7\cdot 13$ &$2\cdot 3\cdot7\cdot 13$&$2\cdot 3\cdot 5\cdot7\cdot 13\cdot 37$&$2\cdot 3\cdot 5\cdot7\cdot 13\cdot 37$\\

\hline $x^2+4$&$17$&$2\cdot5\cdot 17$&$2\cdot 3\cdot5^2\cdot 17$ &$2^2\cdot 3\cdot5^2\cdot 17$&$2^2\cdot 3\cdot 5^2\cdot17\cdot 41$&$2^2\cdot 3\cdot 5^2\cdot13\cdot 17\cdot 41$\\

\hline $x^2+5$&$3\cdot 7$&$2\cdot3\cdot 7$&$2\cdot 3\cdot7\cdot 29$ &$2\cdot 3^2\cdot7\cdot 29$&$2\cdot 3^2\cdot 5\cdot7\cdot 29$&$2\cdot 3^2\cdot 5\cdot7\cdot 29$\\

\hline $x^2+6$&$5^2$&$2\cdot5^2\cdot 7$&$2\cdot 3\cdot 5^2\cdot7\cdot 11$ &$2\cdot 3\cdot 5^2\cdot7\cdot 11$&$2\cdot 3\cdot 5^2\cdot7^2\cdot 11$&$2\cdot 3\cdot 5^2\cdot7^2\cdot 11$\\
\hline
\end{tabular}\end{lrbox}
\resizebox{4.8 in}{!}{\usebox{\tablebox}}

\vskip 1cm

 {\bf Example 3.} Let $f(x)=x^3+b.$ For $1\leq i\leq k,$
we have
$$a_i(x)(x^3+b)+b_i(x)((x+i)^3+b)=C_i$$
where
$$a_i(x)=\left\{\begin{array}{ll}6i^2x^2+(15i^3-9)x+10i^4-18i,&{\rm if}\;3\nmid i\\
6j^2x^2+(45j^3-1)x+90i^4-6i,&{\rm if}\;i=3j.\end{array}\right.$$
$$b_i(x)=\left\{\begin{array}{ll}-6i^2x^2+(3i^3+9)x-i^4-9i,&{\rm if}\;3\nmid i\\
-6j^2x^2+(9j^3+1)x-9j^4-3j,&{\rm if}\;i=3j.\end{array}\right.$$
$$C_i=\left\{\begin{array}{ll}-i^7-27i,&{\rm if}\;3\nmid i\\
-3^5j^7-9j,&{\rm if}\;i=3j.\end{array}\right.$$ Hence, given any
$k\in\mathbb{N}$ and $b\in\mathbb{Z}$, by Theorem \ref{5}, we can
determine the least period $T_{k,f}$ of the arithmetic function
$g_{k,f}.$ For $1\leq k\leq 6$ and $1\leq b\leq 6$, Table II gives
the $T_{k,f}$'s.

\vskip 1cm

\begin{lrbox}{\tablebox}
\begin{tabular}{|c|c|c|c|c|c|c|}
\multicolumn{7}{c}{\bf TABLE II}\\
\multicolumn{7}{c}{$\mbox{ The least period } T_{k,f} \mbox{ of } g_{k,f} \mbox{ with }f(x)=x^3+b$}\\

\multicolumn{7}{c}{}\\

\hline \backslashbox{$f(x)$}{$k$}&1&2&3&4&5&6\\

\hline
 $x^3+1$&$2\cdot7$&$2\cdot
7\cdot13$&$2\cdot3\cdot7\cdot13$&$2^2\cdot3\cdot7\cdot11\cdot13\cdot17\cdot31$&$2^2\cdot3\cdot5\cdot7\cdot13\cdot17\cdot31\cdot43$&$2^2\cdot3\cdot5\cdot7\cdot13\cdot17\cdot19\cdot31\cdot43$\\

\hline
 $x^3+2$&$2\cdot7$&$2\cdot
7\cdot13$&$2\cdot3\cdot7\cdot13$&$2\cdot3\cdot7\cdot11\cdot13\cdot17\cdot31$&$2\cdot3\cdot5\cdot7\cdot13\cdot17\cdot31\cdot43$&$2\cdot3\cdot5\cdot7\cdot13\cdot17\cdot19\cdot31\cdot43$\\

\hline
 $x^3+3$&$2\cdot7$&$2\cdot
7\cdot13$&$2\cdot3\cdot7\cdot13$&$2^2\cdot3\cdot7\cdot11\cdot13\cdot17\cdot31$&$2^2\cdot3\cdot5\cdot7\cdot13\cdot17\cdot31\cdot43$&$2^2\cdot3\cdot5\cdot7\cdot13\cdot17\cdot19\cdot31\cdot43$\\

\hline
 $x^3+4$&$2\cdot7$&$2\cdot
7\cdot13$&$2\cdot3\cdot7\cdot13$&$2\cdot3\cdot7\cdot11\cdot13\cdot17\cdot31$&$2\cdot3\cdot5\cdot7\cdot13\cdot17\cdot31\cdot43$&$2\cdot3\cdot5\cdot7\cdot13\cdot17\cdot19\cdot31\cdot43$\\

\hline
 $x^3+5$&$2\cdot7$&$2\cdot
7\cdot13$&$2\cdot3\cdot7\cdot13$&$2^2\cdot3\cdot7\cdot11\cdot13\cdot17\cdot31$&$2^2\cdot3\cdot5\cdot7\cdot13\cdot17\cdot31\cdot43$&$2^2\cdot3\cdot5\cdot7\cdot13\cdot17\cdot19\cdot31\cdot43$\\

\hline
 $x^3+6$&$2\cdot7$&$2\cdot
7\cdot13$&$2\cdot3\cdot7\cdot13$&$2\cdot3\cdot7\cdot11\cdot13\cdot17\cdot31$&$2\cdot3\cdot5\cdot7\cdot13\cdot17\cdot31\cdot43$&$2\cdot3\cdot5\cdot7\cdot13\cdot17\cdot19\cdot31\cdot43$\\

\hline
\end{tabular}\end{lrbox}
\resizebox{4.8 in}{!}{\usebox{\tablebox}}

$^1$Department of Mathematics, Nanjing University, Nanjing 210093,
P. R. China

{\bf e-mail:} qingzhji@nju.edu.cn

$^2$Department of Mathematics, Nanjing Normal University, Nanjing
210097, P. R. China

{\bf e-mail:} cgji@njnu.edu.cn
\end{document}